\documentclass{article}
\usepackage[utf8]{inputenc}
\usepackage{todonotes}
\usepackage{amsfonts}
\usepackage{amsmath}
\usepackage{amssymb}
\usepackage{mathtools}
\usepackage{graphicx}
\usepackage{amssymb}  
\usepackage{amsthm}   

\usepackage{calc}
\usepackage{hyperref}

\newtheorem{theorem}{Theorem}

\newtheorem{lemma}[theorem]{Lemma}

\newtheorem{remark}[theorem]{Remark}

\def\bu{{\bf u}}

\def\bx{{\bf x}}

\def\basis{{\mathfrak{e}}}

\def\torus{{\mathbb T^2}} 

\title{A pathwise parameterisation for stochastic transport}
\author{Oana Lang \qquad Wei Pan}
\date{}

\begin{document}

\maketitle

\begin{abstract}
   In this work we set the stage for a new probabilistic pathwise approach to effectively calibrate a general class of stochastic nonlinear fluid dynamics models. We focus on a 2D Euler SALT equation, showing that the driving stochastic parameter can be calibrated in an optimal way to match a set of given data. Moreover, we show that this model is robust with respect to the stochastic parameters. 

\end{abstract}

\section{Introduction}
A fundamental challenge in observational sciences, such as weather forecasting and climate change predictions, is
the modelling of uncertainty due, for example, to unknown or neglected physical effects,
and incomplete information in both the data and the formulation of the theoretical models for prediction. 
Various dynamical parameterisation approaches have been proposed to tackle this challenge, see e.g. \cite{Holm2015}, \cite{Wei1}, \cite{Memin2014}, \cite{Wei3}, \cite{EtienneLong}.
Of particular interest are the recently developed Data Driven models, that accommodate uncertainty by predicting both the expected future measurement values and their uncertainties, based on input from measurements and statistical analysis of the initial data.
To effectively incorporate uncertainty in the data driven approach, 
such predictions are made in a probabilistic sense. Additionally, a data assimilation procedure is used to take into account the time integrated information obtained from the data being observed along the solution path during the forecast interval as “in flight corrections”. 

In the geoscience community, \textit{data assimilation} (DA) refers to a set of methodologies designed to efficiently combine past knowledge of a geophysical system (in the form of a \textit{numerical model}) with new information about that system (in the form of \textit{observations}). DA is a central component of Numerical Weather Prediction where it is used to improve forecasting by adjusting the model parameters and reducing the uncertainties. To achieve this, a stochastic feedback loop between the model and the observation may be introduced: the assimilation of more data during the prediction interval will then decrease the uncertainty of the forecasts based on the initial data, by selecting the more likely paths as more observational data is collected. This is
the basis of the so-called \textit{ensemble data assimilation} which uses a set of model trajectories that are intermittently updated according to data. 

A key step for ensuring the successful application of the combined stochastic parameterisation and data assimilation procedure, is the ``correct" calibration of stochastic model parameters.
For Stochastic Advection by Lie Transport (SALT) and Location Uncertainty (LU) models, 
current numerical methods for calibration, see \cite{Wei1}, \cite{EtienneLong}, \cite{Wei3}, have largely been inspired by the physical interpretation of the models derivations. More specifically on the assumption that the flow map is decoupled into a slow scale mean part and a fast scale fluctuating part. In the references mentioned before, it was shown that these methods are effective and led to successful combination of  data driven models and state of the art data assimilation techniques.

In this work, we wish to investigate the feasibility and viability of probabilistic pathwise approach for calibration.
Our general aim is to explore such ideas for a wide class of nonlinear stochastic transport models. This will be very useful in data assimilation problems, as in real world applications the signal is usually observed through discrete observations, but no results of this type for SALT or LU models have been obtained before. Currently, Lagrangian particle trajectories are simulated starting from each point on both the physical grid and its refined version, then the differences between the particle positions are used to calibrate the noise. This is computationally expensive and not fully justified from a theoretical perspective. In the same spirit as \cite{Carsten2} but with a more complicated noise term and without any smoothing effects of a Laplacian, we propose an approach which uses high-frequency in time and low-frequency in space observations of a single path of the solution, to rigorously infer properties of the stochastic parameters. The knowledge of the noise is crucial for determining the behaviour of the solution and for assessing to what degree the solution of the coarse resolution SPDE deviates from the solution of the fine resolution PDE in the model
reduction procedure, so an optimal calibration of the noise parameters is relevant from both a theoretical and an applied perspective. 

    In this work we look at stochastic calibration for the two-dimensional incompressible Euler equation in vorticity form. This stochastic equation models the local rotation of a fluid flow in the presence of spatial uncertainties and it has been derived from fundamental principles in \cite{Holm2015}. This equation is a key ingredient in modelling phenomena in oceanography and in order to ensure that it efficiently encodes the small-scale variablity in the upper part of the ocean, one needs to specify the stochastic parameters based on real observations. 
    One of main issues in parameter estimation using real data is the fact that the model parameters do not map to observations in a unique way (\textit{model identifiability} problem, see e.g. \cite{identifiabilitypaper}).
     For this reason, we believe that a probabilistic approach is much more suitable. 
     
     The 2D Euler equation in the form derived in \cite{Holm2015} and studied in \cite{Wei1}, \cite{Wei3} and \cite{OanaDan1} reads:
    \begin{equation}\label{initialmaineqn}
d{\omega }_{t}+u_{t}\cdot \nabla {\omega }_{t}dt+\displaystyle\sum_{i=1}^{\infty}\xi_i
\cdot \nabla {\omega }_{t}\circ dW_{t}=0 
\end{equation}
where $u=(u^1,u^2)$ is the fluid velocity, $\omega = curl \ u = \partial_2u_1 - \partial_1u_2$ is the vorticity, $(\xi_i)_i$ are divergence-free time-independent vector fields such that 
\begin{equation} \label{xiassumpt1}
\displaystyle\sum_{i=1}^{\infty} \|\xi_i\|_{k+1,\infty}^2 < \infty
\end{equation}
and $(W^{i})_{i \in \mathbb{N}}$ is a sequence of independent Brownian motions.
Global well-posedness for equation \eqref{initialmaineqn} has been studied in \cite{OanaDan1} and the numerical and data assimilation perspective has been studied in \cite{Wei1} and \cite{Wei3}. In \cite{OanaDan1} the authors have proven that equation \eqref{initialmaineqn} admits a uniques pathwise solution which lives in the Sobolev space $\mathcal{W}^{k,2}(\mathbb{T}^2) \ (k\geq 2)$ when $\omega_0 \in \mathcal{W}^{k,2}(\mathbb{T}^2)$ and can be extended to $L^{\infty}(\mathbb{T}^2)$ when $\omega_0 \in L^{\infty}(\mathbb{T}^2)$. \\ 
 
In this paper we consider the following SPDE on the two-dimensional torus $\mathbb{T}^{2}=\mathbb{R}^{2}/\mathbb{Z}^{2}$, driven by a 1-dimensional Brownian motion $W$:
\begin{equation}\label{maineqn}
d{\omega }_{t}+u_{t}\cdot \nabla {\omega }_{t}dt+\xi
\cdot \nabla {\omega }_{t}\circ dW_{t}=0
\end{equation}%
where $u$ and $\omega$ are as above and $\circ$ denotes Stratonovich integration. We impose the following condition on the stochastic parameter $\xi$, in the same spirit as \eqref{xiassumpt1}:
\begin{equation} \label{xiassumpt2}
 \|\xi\|_{k+2,\infty}^2 < \infty
\end{equation}
with $k\geq 2$. This condition ensures that 
for any $f \in \mathcal{W}^{2,2}(\mathbb{T}^2) \cap \mathcal{W}^{2,\infty}(\mathbb{T}^2)$,
\begin{equation}\label{xiassumptl2}
    \|\xi\cdot\nabla f\|_{2}^2 \leq C\|f\|_{1,2}^2  \ \ \ \ \ \ \  \|\xi\cdot\nabla(\xi\cdot\nabla f)\|_2^2 \leq C\|f\|_{2,2}^2
\end{equation}
\begin{equation}\label{xiassumptlinfty}
    \|\xi\cdot\nabla f\|_{\infty}^2 \leq C\|f\|_{1,\infty}^2  \ \ \ \ \ \ \  \|\xi\cdot\nabla(\xi\cdot\nabla f)\|_{\infty}^2 \leq C\|f\|_{2,\infty}^2.
\end{equation}
\begin{remark}
We can view the stochastic part as a space-time noise $(\xi, W)$ where the spatial component is given by $\xi$ and the time component is a standard Brownian motion. This perspective is many times useful in numerical applications where $(\xi \circ dW_t) \cdot \nabla$ is implemented as a random operator applied to the solution $\omega$. 
\end{remark}

 The problem of parameter estimation, known also as \textit{statistical inference}, is technically challenging for such (infinite-dimensional) SPDEs driven by transport noise, as most methods used in the literature benefit from a diagonalizable structure of the underlying space-covariance matrices. This structure is specific for additive noise and therefore it does not apply in our case. Also, most results are obtained for stochastic variations of the heat equation, which contain a smoothing Laplace operator (see for instance \cite{Carsten2}). Our model does not contain a Laplacian a priori, and therefore we cannot exploit the properties of a heat kernel. These makes the analysis much harder.  \\
 
 \noindent\textbf{Contributions of the paper}\\
 In this work, we focus on equation \eqref{maineqn} from two perspectives: 
 \begin{itemize}
     \item First, we show that the driving stochastic parameter $\xi$ can be calibrated in an optimal way to match a set of high-frequency in time given data. This is done using a forced and damped version of the equation and a parametric form of the stream function and the corresponding stochastic parameter which is implemented using an orthonormal basis. 
     Our technique can be explicitly applied to calibrate the 2D Euler model using \textbf{real} oceanic data and we intend to do this in coming work. 
     \item Second, we show that the original 2D Euler model is \textit{robust} with respect to the stochastic parameters $\xi$ in the sense that if we consider two couples $(\omega^1, \xi^1)$ and $(\omega^2, \xi^2)$ which solve equation \eqref{maineqn}, then the $L^2$ distance between $\omega^1$ and $\omega^2$ can be controlled using the initial conditions and the difference between $\xi^1$ and $\xi^2$ only (see Section \ref{section:robustness}). This is important in applications as it shows that if we consider approximate values for $\xi$, the corresponding model solution remains close to the true solution. 
 \end{itemize}

 \noindent\textbf{Structure of the paper}\\
 In Section \ref{sect:problemformulation} below we present the problem formulation. In Section \ref{sect:methodology} we introduce the methodology. In Section \ref{section:robustness} we prove the robustness of the original model and in Section \ref{sect:numericalresults} we present the numerical results. 

\section{Problem formulation}\label{sect:problemformulation}
 Let $(\Omega, \mathcal{F}, (\mathcal{F}_t)_{t\geq 0}, \mathbb{P})$ be a filtered probability space and $W$ a one-dimensional Brownian motion adapted to the complete and right-continuous filtration $(\mathcal{F}_t)_{t\geq 0}$.

We assume we are given a finite sequence of \emph{high frequency} in time vorticity fields, that are denoted by $\omega^*_{t_i}(x)$, $i=1,\dots, N$, and are adapted to $(\mathcal{F}_t)_{t\geq 0}$. We take the view that the $\omega^*_{t_i}$'s are the given observation \emph{data}. We further assume that $\omega^*_{t_i} \in \mathcal{W}^{k,2}(\torus)$.

Writing $\omega_\xi$ to denote solutions to the model \eqref{maineqn} for a given vector field $\xi$, the generic problem we are interested in is to find a $\xi$ so that solutions to \eqref{maineqn} matches the data as best as possible, i.e.
\begin{equation}\label{eq: generic problem}
    \displaystyle\arg\min_{\xi}\|\omega^{*} - \omega_\xi\|
\end{equation}
for some suitable norm.

The dimension of the observations currently coincides with the number of sources of noise, that is we have a \textit{determined} system. However, in practice this is not always a realistic assumption and in future work we will look at \textit{underdetermined} or \textit{overcomplete} systems i.e. when the number of noise sources is larger than the dimension of the observation operator.

In general, the infinite dimensional optimisation problem \eqref{eq: generic problem} may be too hard to solve in practice. We thus make concrete the form of $\xi$. 
Let $(\basis_{j})_{j\in\mathbb{N}}$ be an orthonormal basis in $L^2(\mathbb{T}^2)$. We assume the following parametric form for the stream function of $\xi$, which is henceforth denoted by $\zeta$,
\begin{equation}
    \zeta(x) = \sum_{j=1}^\infty \alpha_j \basis_j,
\end{equation}
where $\alpha_j$ are reals.
Then
\begin{equation}\label{eq: parametric form}
    \xi(x) = \nabla^\perp \zeta(x) = \sum_{j=1}^\infty \alpha_j \nabla^\perp \basis_j(x)
\end{equation}
and the optimisation problem \eqref{eq: generic problem} then reduces to finding the coefficients $\alpha_j$.    

\section{Methodology}\label{sect:methodology}

We will first introduce a couple of known results for vorticity equations, for further discussions on the topic see \cite{MajdaBertozzi} or \cite{MarchioroPulvirenti}. 
 The link between the vorticity and the velocity vector field in equation \eqref{maineqn} is uniquely established using the Biot-Savart operator $K$:
\begin{equation}\label{biotsavart}
\bu(x) = (K \star \omega)(x) = \displaystyle\int_{\mathbb{T}^2}
K(x-y)\omega(y)dy
\end{equation}
with 
\begin{equation}
K(x) = \nabla^{\perp} G(x) = \displaystyle\sum_{k \in \mathbb{Z}^2\setminus
\{0\}} \frac{ik^{\perp}}{\|k\|^2}e^{ik \cdot x}
\end{equation} 
where $G$ is the Green function of the operator $-\Delta$ on $\mathbb{T}^2$
\begin{equation*}
G(x) = \displaystyle\sum_{k \in \mathbb{Z}^2\setminus \{0\}} \frac{e^{ik
\cdot x}}{\|k\|^2}
\end{equation*}
and $k=(k_1, k_2)$, $k^{\perp} = (k_2, -k_1)$. 
It is known that, for any $k\ge 0$, there exists a constant $C_{k,2}$, independent of $u$ such that 
\begin{equation*} \label{biotsavarteqn}
\|\bu\|_{k+1, 2} \leq C_{k,2}\|\omega\|_{k,2}.
\end{equation*}
 If $\psi : \mathbb{T}^2 \times
[0, \infty) \rightarrow \mathbb{R}$ is a solution for $\Delta\psi = -\omega$
then $\bu = \nabla^{\perp}\psi$ solves $\omega = \rm{curl} \ \bu$, so $\bu =
-\nabla^{\perp}\Delta^{-1} \omega$. The reconstruction of $\bu$ from $\omega$ is ensured by the
incompressibility condition $\nabla \cdot \bu = 0$ and  a periodic,
distributional solution of $\Delta \psi = -\omega$ is given by
\begin{equation*}
\psi(x) = (G \star \omega)(x). 
\end{equation*}
From \eqref{maineqn} we have
\begin{equation}\label{eq: maineqn with parametric xi}
    \omega_t(x) = \omega_0(x) - \displaystyle\int_0^t B_s(x;\omega)
    \ ds - \displaystyle\int_0^t \displaystyle\sum_j \alpha_j \nabla^\perp \basis_j(x) \cdot \nabla \omega_t(x)dW_t
\end{equation}
in which for simplicity of notation, we defined $B_s(x;\omega) :=  \bu_s(x) \cdot \nabla \omega_s(x) $.
Combining \eqref{eq: maineqn with parametric xi} with the Biot-Savart law \eqref{biotsavart} we obtain 
\begin{align}
    \bu_t(x) &= \bu_0(x) - \int_0^t \int_\torus K(x-y) B_s(y;\omega) dy ds 
    - \int_0^t \int_\torus K(x-y) \xi(y)\cdot\nabla \omega_s(y) dy \circ dW_s
\end{align}
Now consider the ``kinetic energy"
\begin{equation}
    e_t := \frac12 \int_\torus |\bu_t|^2 dx
\end{equation}
which is not conserved in the SALT case.
Using It\^o's lemma, we obtain

\begin{equation}
    \begin{aligned}
    e_t - e_0 
    & = - \int_0^t  
    \langle \bu_s, K \star (B_s - \frac{1}{2}\xi \cdot \nabla(\xi\cdot\nabla \omega_s)) \rangle  ds \\
     & \qquad -  \int_0^t \langle \bu_s, K\star (\xi \cdot \nabla \omega_s)\rangle\  dW_s
\end{aligned}
\end{equation}
where $\langle .,. \rangle$ is the standard $L^2(\torus)$ pairing.

For a stochastic process $X_t$ defined on a filtered probability space, its \emph{quadratic variation} is defined by \begin{equation}
    \label{eq: quadratic variation}
    [X]_t := \lim_{\max_j\Delta t_j\rightarrow 0}\sum_{i=1}^n |X_{t_i} - X_{t_{i-1}}|^2, 
\end{equation}
where $t_0 = 0 < t_1 < \dots < t_n = t$ is a partition of the interval $[0,t]$,  $\Delta t_i := |t_i - t_{i-1}|$, and the convergence holds in probability (see e.g. \cite{KaratzasShreve}). In our case, for the semimartingale $e_t$, we have
\begin{equation}
\begin{aligned}
    [e]_t &= \int_0^t \langle \bu_s, K\star (\xi \cdot \nabla \omega_s)\rangle^2 \ ds.
\end{aligned}
\end{equation}
Substituting in the parametric form for $\xi$,
we obtain the following quadratic form
\begin{equation}\label{eq: quadratic form}
\begin{aligned}
    [e]_t 
    &= \sum_{i,j=1}^\infty  \alpha_i \alpha_j \int_0^t      \langle \bu_s,   K\star  (\nabla^\perp \basis_j \cdot \nabla \omega_s)\rangle  \langle \bu_s,   K\star  (\nabla^\perp \basis_i \cdot \nabla \omega_s)\rangle \ ds.
\end{aligned}
\end{equation}

Due to global existence and uniqueness of solutions to \eqref{maineqn}, $[e]_t$ exists globally $\mathbb{P}$-almost surely. Thus the right hand side of \eqref{eq: quadratic form} can be arbitrarily well approximated by its truncation for all $t$ i.e. for a given $\epsilon > 0$, there exists $M_\epsilon$ such that
\begin{equation}
    \left | [e]_t - \sum_{i,j=1}^{M_\epsilon} \alpha_i\alpha_j \int_0^t      \langle \bu_s,   K\star  (\nabla^\perp \basis_j \cdot \nabla \omega_s)\rangle  \langle \bu_s,   K\star  (\nabla^\perp \basis_i \cdot \nabla \omega_s)\rangle \ ds \right | < \epsilon.
\end{equation}
Additionally, from the computational perspective, for any fixed $M_\epsilon$, the linear map \begin{equation}\label{eq: linear map}
    \mathbf{A}_{ij} := \int_0^t      \langle \bu_s,   K\star  (\nabla^\perp \basis_j \cdot \nabla \omega_s)\rangle  \langle \bu_s,   K\star  (\nabla^\perp \basis_i \cdot \nabla \omega_s)\rangle \ ds
\end{equation} that defines the truncated quadratic form is symmetric and positive definite\footnote{$[e]_t$ is strictly positive.}, and thus can be diagonalised by a unitary linear map. Doing so, we obtain the following linear problem
\begin{equation}\label{eq: generic linear system}
[e]_t = \sum_{j=1}^{M_\epsilon} \tilde{\alpha}_j^2 \lambda_j + \epsilon',
\end{equation}
where $\epsilon'$ denotes the truncation error of \eqref{eq: quadratic form}, $\lambda_j$ are the eigenvalues of the associated linear map, and $\tilde{\alpha}_j$'s are the original $\alpha$ values which get rescaled by the unitary matrix from the diagonalisation.

We can estimate $[e]_t$ using the high frequency in time data $\omega^*$ and \eqref{eq: quadratic variation}, assuming the discrete sum converges fast enough,
\begin{equation}
    [e]_t \approx \hat{[e]}_{t,N} := \frac14\sum_{i=1}^N (\int_\torus |\nabla^\perp \Delta^{-1}\omega^*_{t_i}|^2 - |\nabla^\perp \Delta^{-1}\omega^*_{t_{i-1}}|^2 dx  )^2.
\end{equation}
The estimate $\hat{[e]}_{t,N}$ could then be used in \eqref{eq: generic linear system} to get an estimate for the $\tilde{\alpha}$. One could then recover the original $\alpha$'s by applying the unitary linear map that's associated with the diagonalisation of $\mathbf{A}_{ij}$.

\begin{remark}\label{rem: qv for vorticity}
The calculations shown in this section could also be directly applied to the vorticity equation \eqref{eq: maineqn with parametric xi} with little modification. In the numerics section of this work, this is what we did. The linear system for estimation, however, would  would also depend on the space variable. 
Recall that the system \eqref{maineqn} conserves spatial integrals of $\omega_t$ (it is a Casimir of the system, see \cite{Wei1}), the quadratic variation of the spatially integrated $\omega$ would be trivial.
\end{remark}

\section{Robustness}\label{section:robustness}
\begin{theorem}\label{thm:robustness}
Let $\omega^1, \omega^2$ be two solutions of the 2D Euler equation \eqref{maineqn} and $\xi^1, \xi^2$ the corresponding stochastic parameters for each of these two solutions. More precisely, $(\omega^{\ell},\xi^{\ell})$ for $\ell=1,2$ solves
\begin{equation}
d{\omega}_{t}^{\ell}+u_{t}^{\ell}\cdot \nabla {\omega }_{t}^{\ell}dt+\xi^{\ell}
\cdot \nabla {\omega }_{t}^{\ell} dW_{t}=\frac{1}{2}\xi^{\ell} \cdot \nabla\left( \xi^{\ell} \cdot \nabla\omega^{\ell}\right). 
\end{equation}%
Then for any $p\geq 2$ there exist non-negative constants $C=C(p,T)$, $C_{1,p}, C_{2,p}$, such that
\begin{equation}\label{robustnessproperty}
    \mathbb{E}\left[ e^{-\gamma(T)}\displaystyle\sup_{t\in[0,T]}\|\omega_t^1-\omega_t^2\|_2^p\right] \leq C \left ( \|\omega_0^1-\omega_0^2\|_2^p + \|\xi^1-\xi^2\|_2^p \right)
\end{equation}
where
\begin{equation*}
    \gamma(T):= C_{1,p}\displaystyle\int_0^T\|\omega_s^1\|_{k,2}^pds + C_{2,p}T^p.
\end{equation*}
\end{theorem}

\vspace{2mm}
\begin{proof}[\textbf{Proof of Theorem \ref{thm:robustness}}.]
Let $\bar{\omega}:= \omega^1 - \omega^2, \bar{u} = u^1-u^2, \bar{\xi} =\xi^1 - \xi^2$. Then $\bar{\omega}$ satisfies
\begin{equation*}
    d{\bar{\omega}}_t + (\bar{u}_t\cdot\nabla\omega_t^1 + u_t^2 \cdot \nabla {\bar{\omega}_t}) dt + \left( \xi^1\cdot\nabla\omega_t^1 - \xi^2\cdot\nabla\omega_t^2 \right) dW_t=\frac{1}{2}\left( \xi^1\cdot\nabla(\xi^1\cdot\nabla\omega_t^1) - \xi^2\cdot\nabla(\xi^2\cdot\nabla\omega_t^2) \right)dt.
\end{equation*}
By the It\^{o} formula: 
\begin{equation*}\label{itodiff}
\begin{aligned}
d\|\bar{\omega}_t\|_2^2=&-2\langle
\bar{\omega}_t,  \xi^1\cdot\nabla\omega_t^1 - \xi^2\cdot\nabla\omega_t^2\rangle dW_t-2\langle
\bar{\omega}_t, \bar{u}_t \cdot \nabla \omega_t^1 + u_t^2 \cdot \nabla {\bar{\omega}_t}\rangle dt \\ 
&+
\left(\langle \bar{\omega}_t, \xi^1\cdot\nabla(\xi^1\cdot\nabla\omega_t^1) - \xi^2\cdot\nabla(\xi^2\cdot\nabla\omega_t^2) \rangle  + \langle \xi^1\cdot\nabla\omega_t^1 - \xi^2\cdot\nabla\omega_t^2, \xi^1\cdot\nabla\omega_t^1 - \xi^2\cdot\nabla\omega_t^2 \rangle\right)
dt. \end{aligned}
\end{equation*}
The difference of the nonlinear terms is analysed explicitly in \cite{OanaDan1} pp. 9:
\begin{equation*}
 \langle
\bar{\omega}_t, \bar{u}_t \cdot \nabla\omega_t^1 \rangle \leq
\|\bar{\omega}_t\|_2\|\bar{u}_t\|_{4}\|\nabla\omega_t^1\|_{4} \leq
C\|\bar{\omega}_t\|_2^2\|\omega_t^1\|_{k,2}.
\end{equation*}
We used here that $\|\nabla\omega_{t}^1\|_{4} \leq C \|\omega_t^1\|_{k,2}$ and
$\|\bar{u}_t\|_{4} \leq C\|\bar{u}_t\|_{1,2} \leq C \|\bar{\omega}_t\|_{2}$. 
Also, since $u^2$ is divergence-free,  
$\langle
\bar{\omega}_t,  u_t^2 \cdot \nabla {\bar{\omega}_t}\rangle =-\frac{1}{2} \displaystyle\int_{\mathbb{T}^2} (\nabla \cdot u_t^2)
(\bar{\omega}_t)^2 dx =0 $. 
We estimate the difference terms which include $\xi^1$ and $\xi^2$ in Lemma \ref{estimates} below. Note here that the term $\langle \bar{\omega}_t, \xi^2 \cdot \nabla \left( \xi^2\cdot\nabla \bar{\omega}_t\right)\rangle$ is negative. 
To estimate the stochastic term let
\begin{equation*}
    D_t : = \displaystyle\int_0^t\langle
\bar{\omega}_s,  \xi^1\cdot\nabla\omega_s^1 - \xi^2\cdot\nabla\omega_s^2\rangle dW_s
\end{equation*}
and
\begin{equation*}
    m_t := \|\bar{\omega}_t\|_2^2 \ \ \ \ \ \ \ \ \ \ \ \ \ \ \ \ \ \ Z: = \|\bar{\xi}\|_2^2.
\end{equation*}
By the Burkholder-Davis-Gundy inequality, for arbitrary $p\geq 2$ there exists a constant $C_p$ such that\footnote{In this proof $C, C_p$ are generic constants which may differ from line to line and from term to term.}
\begin{equation*}
    \mathbb{E}\left[ \displaystyle\sup_{s\in[0,t]}|D_s|^p\right] \leq C_p \mathbb{E}\left[ [D]_t^{p/2} \right]
\end{equation*}
where $[D]_t$ is the quadratic variation of the martingale $D_t$.
Note that $[D]_t$ can be controlled using the control for $Q$ in Lemma \ref{estimates} so we have
\begin{equation*}
    \begin{aligned}
    \mathbb{E}\left[ [D]_t^{p/2}\right] \leq C_p \displaystyle\int_0^t \mathbb{E}\left[ \displaystyle\sup_{r\in[0,s]}\left(m_r^p + \tilde{Z}^p\right)\right] ds
    \end{aligned}
\end{equation*}
where $\tilde{Z}:= C\|\omega_t^1\|_{k,2}^2Z$ with $k\geq 3$.
We then have, for $t \in [0,T]$:
\begin{equation*}
    m_t \leq m_0 + \displaystyle\int_0^t\psi(s)ds -2D_t +  \displaystyle\int_0^t \phi(s)m_s ds
\end{equation*}
where
\begin{equation*}
   \phi(t) := C\|\omega_t^1\|_{k,2}^2 + \tilde{C} 
\end{equation*}
and
\begin{equation*}
    \psi(t) := (C\|\omega_t^1\|_{k,2}^2 + 1)Z.
\end{equation*}
Then 
\begin{equation*}
    \begin{aligned}
    d\left( e^{-\displaystyle\int_0^t\phi(r)^pdr}m_t^p\right) & = e^{-\displaystyle\int_0^t\phi(r)^pdr}\left( dm_t^p - m_t^p\phi(t)^pdt\right) \\
    & \leq e^{-\displaystyle\int_0^t\phi(r)^pdr} \left( \psi(t)^pdt - 2(dD_t)^p\right)
    \end{aligned}
\end{equation*}
and
\begin{equation*}
  \begin{aligned}
  \mathbb{E}\left[ e^{-\displaystyle\int_0^t\phi(r)^pdr} \displaystyle\sup_{s\in[0,t]}m_s^p \right] &\leq C_p\left( m_0^p + \mathbb{E}\left[\displaystyle\int_0^te^{-\displaystyle\int_0^s\phi(r)^pdr}\displaystyle\sup_{r\in[0,s]}\psi(r)^pds\right] + 2\mathbb{E}\left[e^{-\displaystyle\int_0^t\phi(r)^pdr}\displaystyle\sup_{s\in[0,t]}|D_s|^p\right] \right)\\
  & \leq  C_p\left(m_0^p + \mathbb{E}\left[\displaystyle\int_0^te^{-\displaystyle\int_0^s\phi(r)^pdr}\displaystyle\sup_{r\in[0,s]}\psi(r)^pds\right] \right)\\
  & \leq C_p\left(m_0^p + \mathbb{E}\left[\displaystyle\int_0^t\displaystyle\sup_{r\in[0,s]}\psi(r)^pds\right] \right)
\end{aligned}
\end{equation*}
since $D_t$ is a martingale so its expectation is zero. 
Therefore
\begin{equation*}
    \begin{aligned}
    \mathbb{E}\left[e^{-\gamma(t)}\displaystyle\sup_{s\in[0,t]}m_s^p \right] \leq 
    C_p\left( m_0^p +  \mathbb{E}\left[\displaystyle\int_0^t\displaystyle\sup_{r\in[0,s]}\psi(r)^pds\right]\right) 
    \end{aligned}
\end{equation*}
where $\gamma(t):= \displaystyle\int_0^t \phi(s)^p ds$.
This gives
\begin{equation*}
    \mathbb{E}\left[e^{-\gamma(T)}\displaystyle\sup_{s\in[0,T]}\|\omega_s^1-\omega_s^2\|_2^p\right] \leq C_{p,T}\left( \|\omega_0^1-\omega_0^2\|_0^p + \|\xi^1-\xi^2\|_2^p\right)
\end{equation*}
for any $p\geq 2$, with $\gamma$ as above.
We used here the fact that the 2D Euler equation \eqref{maineqn} has a unique global solution in $\mathcal{W}^{k,2}(\mathbb{T}^2)$ for $k\geq 2$, as proven in \cite{OanaDan1}. Moreover $\bar{\xi}$ is deterministic and time-independent so 
\begin{equation*}
\begin{aligned}
\mathbb{E}\left[\displaystyle\int_0^T \psi(s)ds\right]&=\mathbb{E}\left[ \displaystyle\int_0^T (C\|\omega_s^1\|_{k,2}^2 +1)\|\bar{\xi}\|_2^2ds\right] \\
& \leq \|\bar{\xi}\|_2^2 \left(CT\mathbb{E}\left[ \displaystyle\sup_{s\in[0,T]}\|\omega_s^1\|_{k,2}^2\right] + T \right) \\
& \leq \tilde{C}_T \|\bar{\xi}\|_2^2.
\end{aligned}
\end{equation*}
\end{proof}
\begin{remark}
Note that we have shown above a result which is similar to an application of the Growall Lemma \ref{gronwall} in \ref{appendix}. However, since we work with random variables, we made this step explicitly, for clarity. 
\end{remark}

\begin{lemma}\label{estimates}
Let $(\omega_t^1,\xi^1)$ and $(\omega_t^2,\xi^2)$ be two solutions of the 2D Euler equation with $\bar{\omega}_t:= \omega_t^1 - \omega_t^2$ and $\bar{\xi}:= \xi^1 -\xi^2.$ Then there exists constants $C$\footnote{$C$ differs from line to line and from term to term depending on the Sobolev embedding we use.} such that the following estimates hold:
\end{lemma}

\begin{subequations}
\vspace{-3mm}
\begin{equation}
    Q := |\langle
\bar{\omega}_t,  \xi^1\cdot\nabla\omega_t^1 - \xi^2\cdot\nabla\omega_t^2\rangle| \leq C\|\bar{\omega}_t\|_2^2 + C\|\omega_t^1\|_{k,2}^2\|\bar{\xi}\|_{2}^2.
\end{equation}
\begin{equation}
    A :=  \langle \xi^1 \cdot \nabla\omega_t^1 - \xi^2 \cdot \nabla \omega_t^2, \xi^1 \cdot \nabla\omega_t^1 - \xi^2 \cdot \nabla \omega_t^2 \rangle \leq  C \|\bar{\omega}_t\|_2^2 + C\|\omega_t^1\|_{k,2}^2 \|\bar{\xi}\|_2^2 
\end{equation}
\end{subequations}
\begin{equation}
    |B| \leq C\|\omega_t^1\|_{k,2}^2\|\bar{\omega}_t\|_2^2 + \|\bar{\xi}\|_2^2 
\end{equation}
where 
\begin{equation*}
    B  :=  \langle \omega_t^1 - \omega_t^2, \xi^1 \cdot \nabla \left(\xi^1\cdot\nabla\omega_t^1 \right) - \xi^2 \cdot \nabla \left(\xi^2\cdot\nabla\omega_t^2 \right) \rangle. 
\end{equation*}
and $k\geq 4$. 
\vspace{2mm}
\begin{proof}
For the difference terms which include $\xi^1$ and $\xi^2$ we use that
\begin{equation*}
   \begin{aligned}
   \xi^1 \cdot \nabla\omega_t^1 - \xi^2 \cdot \nabla \omega_t^2 = \bar{\xi} \cdot \nabla\omega_t^1 + \xi^2 \cdot \nabla \bar{\omega}_t.
   \end{aligned}
\end{equation*}
We have
\begin{equation*}
    \begin{aligned}
    Q &= |\langle
\bar{\omega}_t,  \xi^1\cdot\nabla\omega_t^1 - \xi^2\cdot\nabla\omega_t^2\rangle| \\
     & \leq |\langle \omega_t^1-\omega_t^2, (\xi^1-\xi^2)\cdot\nabla \omega_t^1 \rangle| + |\langle \omega_t^1-\omega_t^2, \xi^2\cdot\nabla(\omega_t^1-\omega_t^2)\rangle|\\
     & \leq \frac{1}{2}\|\omega_t^1-\omega_t^2\|_2^2 + \frac{1}{2}\|\nabla\omega_t^1\|_{\infty}^2\|\xi^1-\xi^2\|_{2}^2 \\
     & \leq \frac{1}{2}\|\bar{\omega}_t\|_2^2 + \frac{C}{2}\|\omega_t^1\|_{k,2}^2\|\bar{\xi}\|_{2}^2
    \end{aligned}
\end{equation*}
with $k\geq 3$, 
since the second scalar product is zero due to the fact that $\div \xi^2 = 0.$
Also 
\begin{equation*}
    \begin{aligned}
  A &=  \langle \xi^1 \cdot \nabla\omega_t^1 - \xi^2 \cdot \nabla \omega_t^2, \xi^1 \cdot \nabla\omega_t^1 - \xi^2 \cdot \nabla \omega_t^2 \rangle = \|\xi^1 \cdot \nabla\omega_t^1 - \xi^2 \cdot \nabla \omega_t^2\|_2^2 \\
  & \leq \|(\xi^1 - \xi^2) \cdot \nabla\omega_t^1\|_2^2 + \|\xi^2 \cdot \nabla (\omega_t^1 - \omega_t^2)\|_2^2 \\
  & \leq \|\xi^1 - \xi^2\|_2^2 \|\nabla\omega_t^1\|_{\infty}^2 + C \|\omega_t^1-\omega_t^2\|_2^2 \\
 & \leq C\|\omega_t^1\|_{k,2}^2\|\bar{\xi}\|_2^2  + C \|\bar{\omega}_t\|_2^2 
\end{aligned}
\end{equation*}
where $k\geq 3$. 
For the higher order term we have
\begin{equation*}
    \begin{aligned}
    B & =  \langle \omega_t^1 - \omega_t^2, \xi^1 \cdot \nabla \left(\xi^1\cdot\nabla\omega_t^1 \right) - \xi^2 \cdot \nabla \left(\xi^2\cdot\nabla\omega_t^2 \right) \rangle  \\
      & =  \langle \omega_t^1 - \omega_t^2, (\xi^1 - \xi^2) \cdot \nabla (\xi^1 \cdot \nabla \omega_t^1) \rangle  \\
      & + \langle \omega_t^1 - \omega_t^2, \xi^2 \cdot \nabla \left( (\xi^1-\xi^2) \cdot \nabla \omega_t^1\right) \rangle \\
      & + \langle \omega_t^1 - \omega_t^2, \xi^2 \cdot \nabla \left( \xi^2\cdot\nabla(\omega_t^1 - \omega_t^2)\right)\rangle \\
      & =: a + b + c. 
    \end{aligned}
\end{equation*}
Note that $c$ is negative: 
\begin{equation*}
    \begin{aligned}
    \langle \omega_t^1 - \omega_t^2, \xi^2 \cdot \nabla \left( \xi^2\cdot\nabla(\omega_t^1 - \omega_t^2)\right)\rangle & = - \langle \xi^2 \cdot \nabla(\omega_t^1 - \omega_t^2),   \xi^2\cdot\nabla(\omega_t^1 - \omega_t^2) \rangle \\
    & = -\|\xi^2 \cdot \nabla(\omega_t^1 - \omega_t^2)\|_2^2 \\
    & \leq 0 
    \end{aligned}
\end{equation*}
so $|B| \leq |a| + |b|$.
We estimate $|a|$ as follows:
\begin{equation*}
    \begin{aligned}
  |\langle \omega_t^1 -\omega_t^2, (\xi^1 -\xi^2)\cdot \nabla( \xi^1 \cdot \nabla \omega_t^1)\rangle|  &\leq \frac{1}{2} \|\nabla(\xi^1\cdot\nabla\omega_t^1)\|_{\infty}^2\|\omega_t^1-\omega_t^2\|_2^2 + \frac{1}{2}\|\xi^1-\xi^2\|_2^2 \\
  & \leq \frac{C}{2} \|\omega_t^1\|_{2,\infty}^2\|\omega_t^1-\omega_t^2\|_2^2 + \frac{1}{2}\|\xi^1-\xi^2\|_2^2 \\
  & \leq \frac{C}{2} \|\omega_t^1\|_{k,2}^2\|\bar{\omega}_t\|_2^2 + \frac{1}{2}\|\bar{\xi}\|_2^2
    \end{aligned}
\end{equation*}
with $k\geq 4$.
Likewise, we estimate $|b|$: 
\begin{equation*}
    \begin{aligned}
    |\langle \omega_t^1 - \omega_t^2, \xi^2 \cdot \nabla \left( (\xi^1-\xi^2) \cdot \nabla \omega_t^1\right) \rangle|  & = |\langle \xi^2\cdot\nabla(\omega_t^1 - \omega_t^2),   (\xi^1-\xi^2) \cdot \nabla \omega_t^1 \rangle| \\ 
    & \leq \frac{1}{2} \|\nabla\omega_t^1\|_{\infty}^2\|\xi^2\cdot\nabla(\omega_t^1-\omega_t^2)\|_2^2 + \frac{1}{2}\|\xi^1-\xi^2\|_2^2 \\
    & \leq \frac{C}{2} \|\nabla\omega_t^1\|_{\infty}^2\|\omega_t^1-\omega_t^2\|_2^2 + \frac{1}{2}\|\xi^1-\xi^2\|_2^2 \\
    & \leq \frac{C}{2} \|\omega_t^1\|_{k,2}^2\|\bar{\omega}_t\|_2^2 + \frac{1}{2}\|\bar{\xi}\|_2^2
    \end{aligned}
\end{equation*}
with $k\geq 3$.
\end{proof}



\section{Numerical results}\label{sect:numericalresults}

We implemented the main equation \eqref{maineqn} with added forcing and damping, on a unit square domain with doubly periodic boundary conditions,
\begin{equation}\label{eq: main numerical eqn}
d{\omega }_{t}+\bu_{t}\cdot \nabla {\omega }_{t}dt+\xi
\cdot \nabla {\omega }_{t}\circ dW_{t}= (Q - r\omega_t) dt 
\end{equation}%
where we chose $r=0.001$ and $Q(x) = 0.01 (\cos(8 \pi y) + \sin( 8 \pi x) )$. 
We considered a $\xi$ whose parametric form with respect to the Fourier basis 
consists of only one $\alpha$. The stream function of our chosen $\xi$ is given by 
\begin{equation}\label{eq: zeta}
\zeta(x, y) = \alpha \left(\cos(k_1 2\pi x)\cos(k_2 2\pi y) - \sin(k_1 2\pi x)\sin(k_2 2\pi y)\right).
\end{equation}
Note that
\begin{equation}
    \zeta = \frac\alpha2(e^{i2\pi k\cdot x } + e^{-i2\pi k\cdot x}),
\end{equation}
and
\begin{equation}
    \xi = i\alpha \pi (e^{i2\pi k\cdot x } - e^{-i2\pi k\cdot x})k^\perp.
\end{equation}

\begin{figure}[ht!]
\centering
\includegraphics[width=0.9\textwidth]{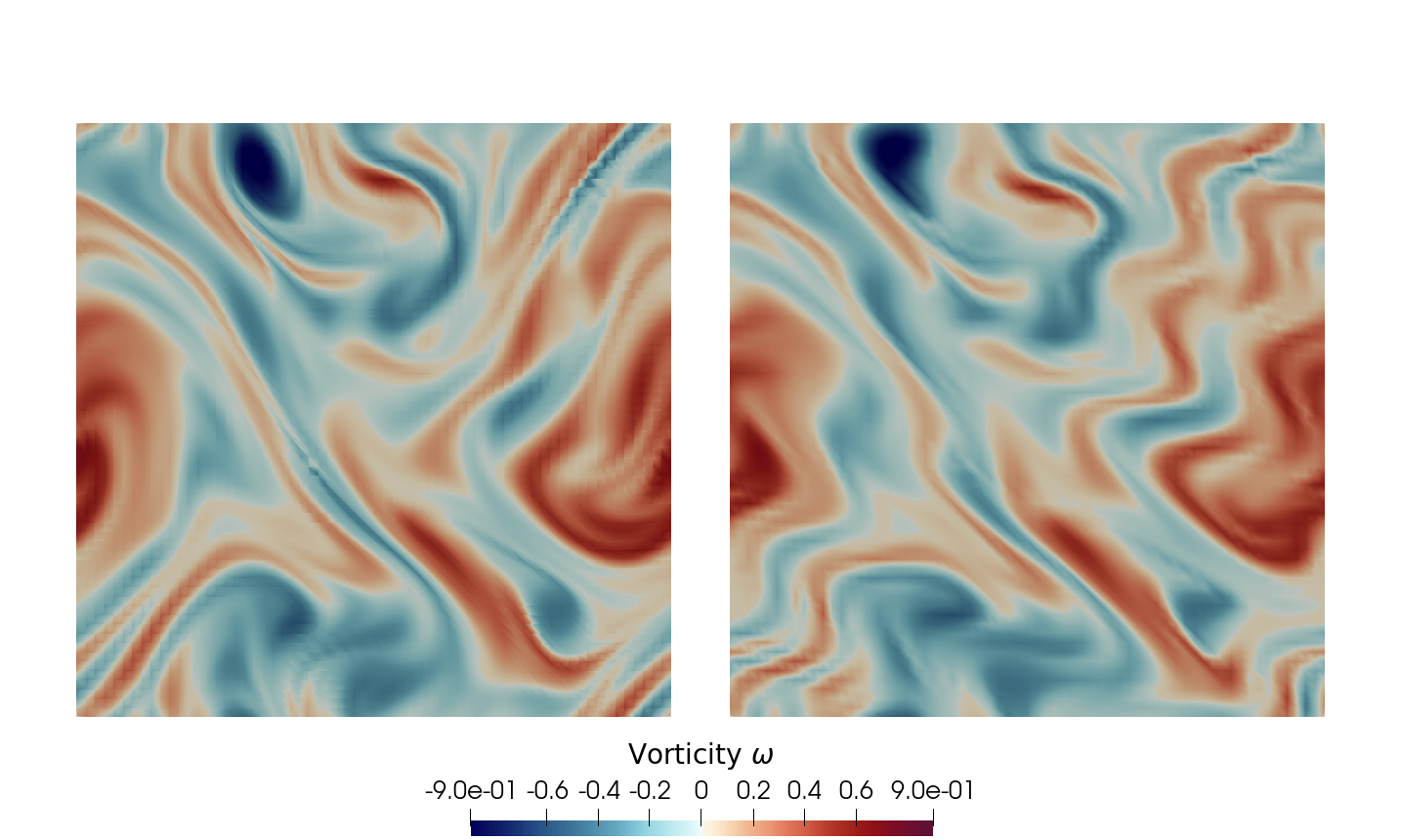}
\caption{Snapshots of the numerical solution $\omega(t,x)$ to \eqref{eq: main numerical eqn} at times $t=0$ (left), and $t=1$ (right). }\label{fig: vorticity snapshots}
\end{figure}

To discretise \eqref{eq: main numerical eqn}, we followed the methods documented in \cite{Wei1} -- a mixed Finite Element method was used for the spatial derivatives, and an explicit strong stability preserving Runge-Kutta scheme of order 3 was used for the time derivative. 
We added the forcing and damping terms to balance out the energy dissipation caused by discretisation. This helped with maintaining the statistical homogeneity of the numerical solution, once it has reached a spun-up state from some set initial state. 
Our choice for the set initial state was
\begin{equation}\label{eq: initial state}
\begin{aligned}
    \omega(0,x,y) &= \sin(8\pi x)\sin(8\pi y) + 0.4 \cos(6\pi x)\cos(6 \pi y) \\&\quad + 0.3\cos(10 \pi x)\cos(4\pi y) + 0.02\sin(2\pi y) + 0.02\sin(2\pi x).
\end{aligned}
\end{equation}
Spatially, we chose a grid of size $64 \times 64$ cells. We first spun-up the system until it reached a statistical equilibrium state. This statistical equilibrium state was then set as the initial condition for our experiment. Figure \ref{fig: vorticity snapshots} shows a snapshot of the obtained initial condition. 
Over the spin-up phase, we used a smaller $\alpha=0.000001$ value and $k^\intercal =(2,4)$.

The time horizon for the experiment data was chosen to be the unit interval, i.e. we generated data $\omega^*(t_i,x)$ for $0=t_0 < t_1 < \dots < t_N = 1$. See Figure \ref{fig: vorticity snapshots} for snapshots of $\omega^*(0, \bx)$ and $\omega^*(1, \bx)$. When generating the data, we used the larger value of
$\alpha=0.001$. This was so to avoid any possible numerical issues\footnote{When $\alpha$ is small, $\alpha^2$ is close to machine precision.} when we attempted to recover $\alpha$ from data.

For these concrete experiment parameter choices, we chose to follow Remark \ref{rem: qv for vorticity}, and directly worked with the quadratic variation of $\omega$.

Assuming we know in-advance the exact Fourier wavenumber $k$, the linear system for estimation reduces to 
\begin{align}
  \hat{[\omega]}_{t,N}(\bx) := \displaystyle\sum_{i=1}^N(\omega_{t_i}(\bx) - \omega_{t_{i-1}}(\bx))^2 
\approx \alpha^2 4 \pi^2 \
B(t,k,\bx) 
\basis'_k(\bx)
\end{align}
where 
\begin{equation}
    B(t,k,\bx):= \displaystyle\int_0^t ( k^\perp \cdot \nabla\omega_s(\bx))^2ds 
\end{equation}
and
\begin{equation}
\basis'_k(x,y) := \left(
\cos(k_1 2\pi x)\sin(k_2 2\pi y)+\sin(k_12\pi x)\cos(k_2 2\pi y)
\right)^2.    
\end{equation}
Thus our estimate for $\alpha$ is given by
\begin{equation}\label{eq: numerical estimate for alpha}
    \hat{\alpha}_N^2 = \frac1{4\pi^2}\frac{\int_\torus \hat{[\omega]}_{t,N}(\bx) d\bx }{\int_\torus B(t,k,\bx)\basis'_k(\bx) d\bx}.
\end{equation}

\begin{remark}
In \eqref{eq: numerical estimate for alpha}, we applied spatial averaging before dividing to avoid possible division by zero issues, and to help with stablising estimation errors.
\end{remark}

\begin{remark}
The assumption that we know $k$ in advance is of course too strong from the applications viewpoint. The aim of this experiment is to test the strength of the pathwise approach under the assumption of ''perfect knowledge". If we cannot accurately recover $\alpha$ in this case, then getting a good estimate for $\alpha$ using the pathwise approach may be too difficult or impractical in more realistic scenarios.
\end{remark}

\begin{figure}[ht!]
\centering
\includegraphics[width=.9\textwidth]{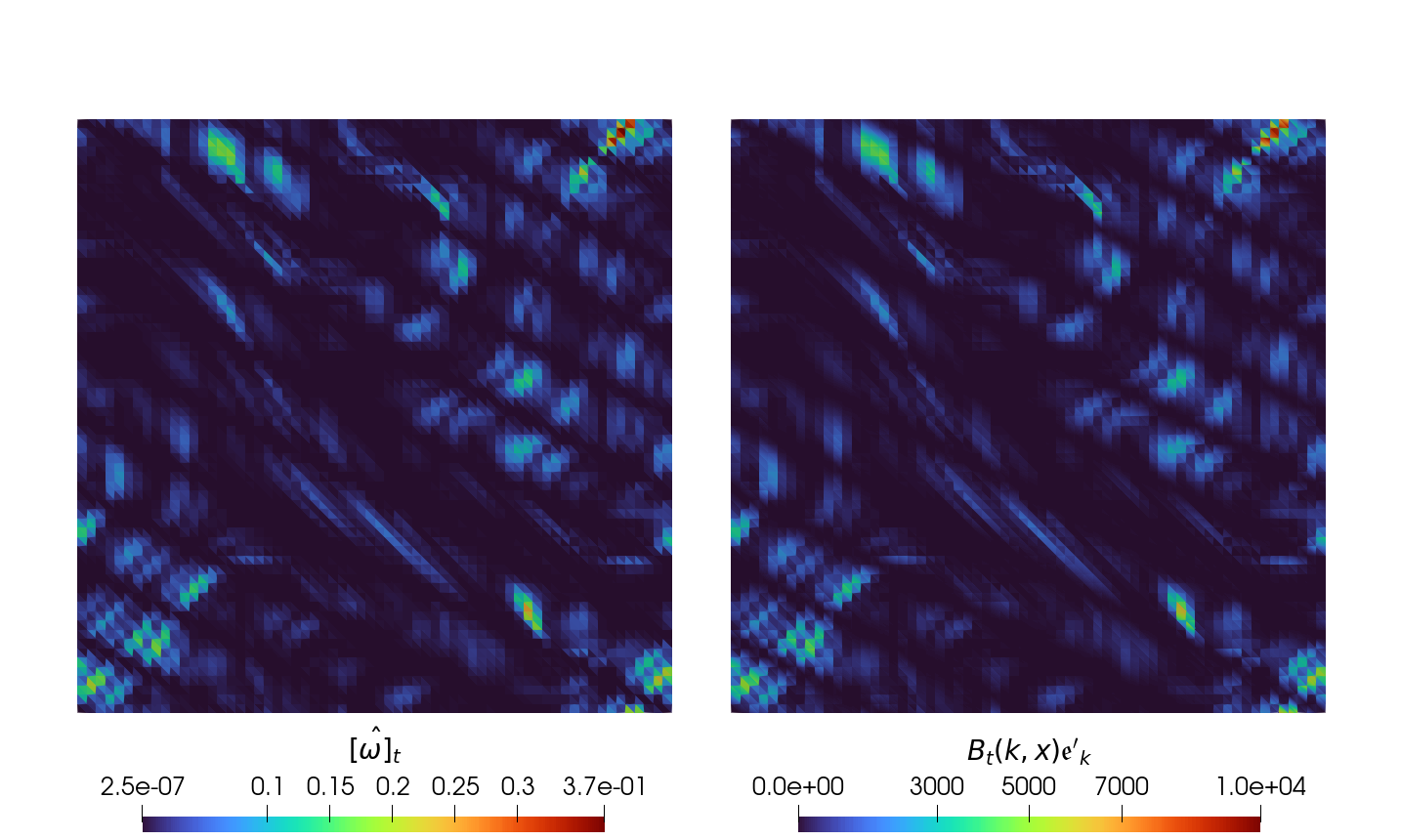}
\caption{Shown on the left is a snapshot of the estimate $\hat{[\omega]}_t$, which was computed using $N=200000$ data samples. Shown on the right is a snapshot of the basis element $B_t(k,x)\left(
\cos(k_1 2\pi x)\sin(k_2 2\pi y)+\sin(k_12\pi x)\cos(k_2 2\pi y)
\right)^2$, which was approximated using the same $N$ number of data samples.} \label{fig: linear system picture}
\end{figure}

\begin{figure}[ht!]
\centering
\includegraphics[width=.8\textwidth]{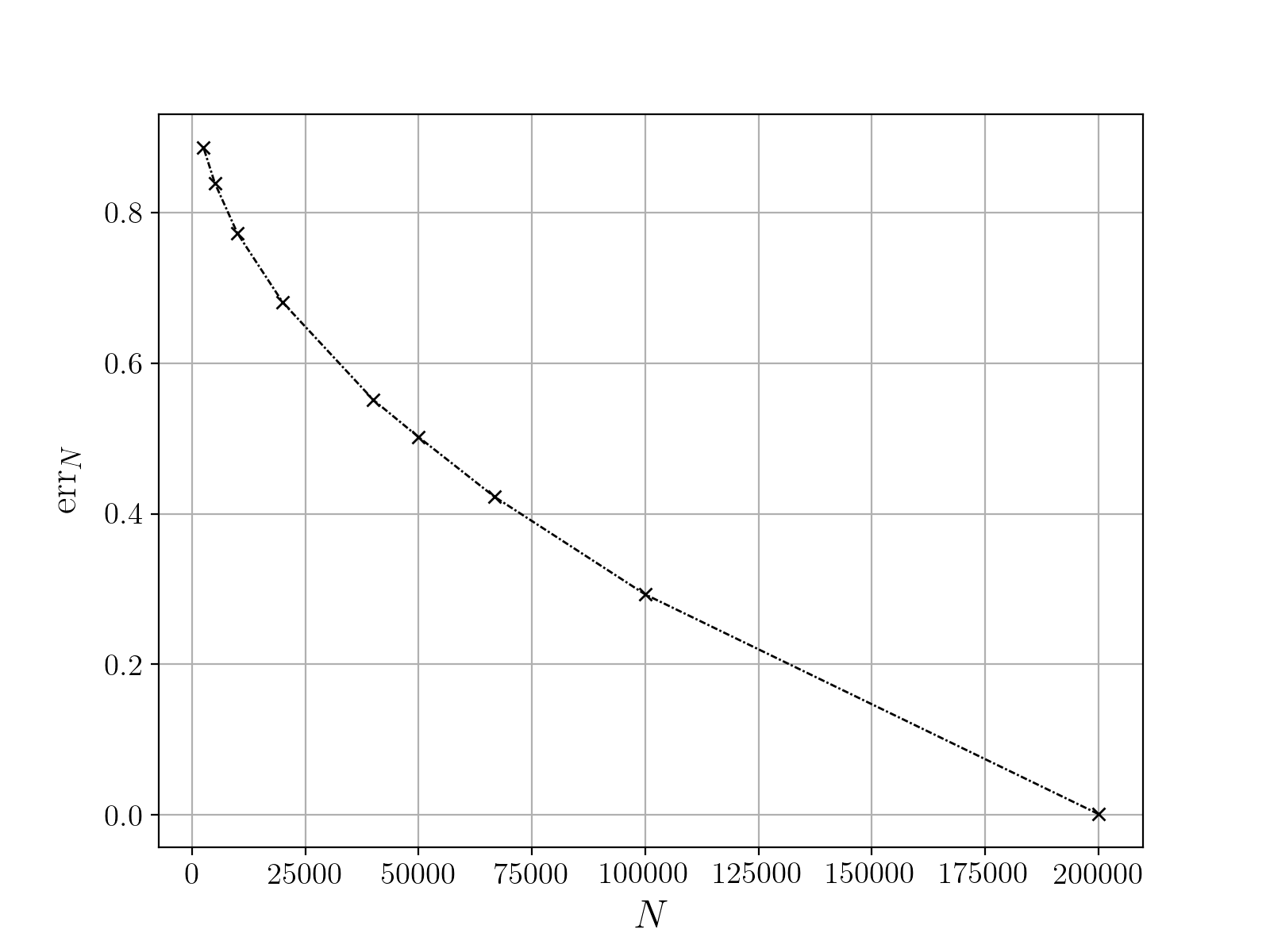}
\caption{The plot shows the relative error $\mathrm{err}_N$ defined in \eqref{eq: relative error} as a function of $N$. $\mathrm{err}_N$ was computed for $N=2500, 5000, 10000, 20000, 40000, 50000, 66667, 100000, 200000.$} \label{fig: relative error}
\end{figure}

Figure \ref{fig: linear system picture} shows snapshots of $\hat{[\omega]}_{t,N}(\bx)$ and $B(t,k,\bx)\basis'_k(\bx)$. We applied  \eqref{eq: numerical estimate for alpha} for different values of $N$. In each case, the time integral that constitutes $B(t,k,\bx)$ was approximated using a simple trapezoidal rule, for which the same $N$ number of data snapshots were used. Figure \ref{fig: relative error} shows the results for the relative error
\begin{equation}\label{eq: relative error}
    \mathrm{err}_N = \frac{|\alpha - \hat{\alpha}_N|}{\alpha}
\end{equation}
for the different values of $N$. The results show that, in the worst case of $N=2500$, the relative error was no greater than 0.89. This translates to an absolute error of range of $0.001 \pm 0.00089$. The best case was when all $200000$ data samples were used to estimate $\alpha$, the relative error in that case was $0.00135$. This suggests convergence and stabilisation of the sum for $\hat{[\omega]}_t$.

For future work, we aim to test the pathwise approach for cases in which we do not know the exact selection of basis elements for $\xi$. Further, we wish to extend and test these ideas on coarse grained PDE data and compare with the results that were obtained in \cite{Wei1} using previously developed calibration methods.

\noindent\textbf{Acknowledgments} \\
The authors would like to thank Prof Dan Crisan for the many helpful suggestions and constructive ideas he shared with them during the preparation of this work. They also thank Prof Darryl Holm, Prof Bertrand Chapron, Prof Etienne Mémin, and the whole STUOD team for many inspiring discussions they had during the STUOD meetings. \\

\noindent\textbf{Funding} \\
Both authors were partially supported by the European Research Council (ERC) under the European Union’s Horizon 2020 Research and Innovation Programme (ERC, Grant Agreement No 856408).

\section{Appendix}\label{appendix}
\begin{lemma}{(Gronwall lemma)}\label{gronwall}
Let $\beta : [0,T] \rightarrow [0,\infty)$ be a non-negative absolutely continuous function that satisfies for a.e. $t$
\begin{equation*}
    d\beta(t) \leq \phi(t)\beta(t)dt + \psi(t)dt
\end{equation*}
where $\phi,\psi$ are non-negative integrable functions on $[0,T]$. Then
\begin{equation*}
    \beta(t) \leq e^{\displaystyle\int_0^t\phi(s)ds}\left( \beta(0) + \displaystyle\int_0^t \psi(s)ds\right)
\end{equation*}
for all $t\in[0,T]$. 
\end{lemma}

\end{document}